\documentclass[11pt]{article}%
\usepackage{amssymb}
\newtheorem{theorem}{Theorem}%
\newtheorem{lemma}[theorem]{Lemma}%
\newtheorem{rem}[theorem]{Remark}%
\def\lg{\langle} \def\rg{\rangle}

\def\f{\noindent}

\def\qed{\hfill $\Box$}

\def\demo{{\bf Proof}\hskip10pt}

\def\N{\hbox{\rm N}}

\def\C{\hbox{\rm C}}

\def\Z{\hbox{\rm Z}}
\def\H{\hbox{\rm H}}

\def\cl{\hbox{\rm cl}}

\begin{document}
\baselineskip=16pt
\title{A refinement on a theorem of Z. Janko
\thanks{This work was supported by NSFC(Nos. 11901367, 11971280, 11771258).}}
\author{\\ Qiangwei Song, Lijian An\thanks{Corresponding author. e-mail: anlj@sxnu.edu.cn}\\
Mathematics Department, Shanxi Normal University\\
Linfen Shanxi 041004 China\\
\\}
\date{}
\maketitle

\begin{abstract}
We say that a subgroup $H$ is isolated in a group $G$ if for each $x\in G$ we have either $x\in H$ or $\langle x\rangle \cap H={1}$.
Z. Janko, in his paper [J. Algebra, 465(2016), 41--61], determined certain classes of finite nonabelian $p$-groups which possess some isolated subgroups.
In this note, a theorem of his paper is refined.

\medskip
\vspace{0.2cm}

 \f {\bf Key Words:} Finite $p$-groups; Isolated subgroups; Minimal nonabelian subgroups; Hughes subgroup

\medskip
\f {\it \bf AMS Subject Classifications:} 20D15.
\end{abstract}

%\section{Introduction}

In this note, we refine two results of Janko, which are \cite[Theorem 4, Corollary 5]{Janko}. The notation are coincide with \cite{Janko}. Now we need the following preliminary lemmas.

%we will prove that $p\geqslant 5$ for the possibility (a) in Theorem 4, the possibility (b) in Theorem 4 does not occur, and the $p$-groups satisfying
%the hypothesis of \cite[Corollary 5]{Janko} do not exist.
%%Moreover, We will give a necessary and sufficient condition.
%In other words, Theorem 4 and Corollary 5 in \cite{Janko} can be refined as follows:

\begin{lemma}\label{p>2}
Let $G$ be a finite $2$-group. If $G$ has a proper subgroup $H$ which is isolated in $G$, then $H$ is abelian.
\end{lemma}
\demo Let $T\leq G$ such that $|T:H|=2$.
Since $H$ is isolated in $G$, $H$ is isolated in $T$ and so $o(x)=2$ for all $x\in T-H$.
Let $y\in H$. Then $o(xy)=2$. So $y^x=y^{-1}$.  Hence $(ab)^x=(ab)^{-1}$ for all $a,b\in H$. On the other hand,
$(ab)^x=a^xb^x=a^{-1}b^{-1}$. So $ab=ba$. i.e., $H$ is abelian.\qed

\begin{lemma}\label{gorder=3}
Suppose that $G$ is a finite $3$-group, $M$ is maximal in $G$ and $g^3=1$ for all $g\in G-M$. If $[a,b]=1$ for some $a,b\in M$, then $[a^x,b]=1$ for all $x\in G-M$.
\end{lemma}
\demo Let $x\in G-M$ and $y\in M$. Then $(xy)^3=1$. That is,  $yxy=x^{-1}y^{-1}x^{-1}$. Now
\begin{eqnarray*}
[a^x,b]&=&x^{-1}a^{-1}xb^{-1}x^{-1}axb\\
&=&x^{-1}a^{-1}(xb^{-1}x)(xax)b\\
&=&x^{-1}a^{-1}(bx^{-1}b)(a^{-1}x^{-1}a^{-1})b\\
&=&x^{-1}a^{-1}bx^{-1}a^{-1}bx^{-1}a^{-1}b\\
&=&(x^{-1}a^{-1}b)^3=1.
\end{eqnarray*}\qed

\begin{lemma}\label{njhqgl}
Suppose that $G$ is a finite $p$-group, and $H$ is a proper minimal nonabelian subgroup of $G$. If $H$ is isolated in $G$, then $H$ is isomorphic to ${\rm S}{\rm (}p^3{\rm )}${\rm(}a nonabelian group of order $p^3$ and exponent $p${\rm)}.
\end{lemma}
\demo  By Lemma \ref{p>2}, we have $p>2$. Let $T\leq G$ such that $|T:H|=p$.
Since $H$ is isolated in $G$, $H$ is isolated in $T$ and so all elements in $T-H$ are of order $p$.

If $p\geq5$, then by \cite[Corollary $3.4$]{QXA}, $T$ has at least two minimal nonabelian subgroups of index $p$. Moreover, by \cite[Theorem $2.7$]{An4}, we have $\cl(G)\leq 3$. It follows from \cite[Theorem $7.1$]{Ber1} that $T$ is regular.
Notice that $T=\lg T-H\rg$ and all elements in $T-H$ are of order $p$. We have $\exp(T)=p$. Moreover, $\exp(H)=p$. Hence $H\cong {\rm S}{\rm (}p^3{\rm )}$.

 In the following, assume $p=3$.

First, we claim that $T'\le \Z(H)=\Phi(H)$. In fact, since $H$ is minimal nonabelian, $\Z(H)=\Phi(H)$ and $|H:\Z(H)|=3^2$ by \cite[Lemma $2.2$]{ZAX}.
Let $h\in H-\Z(H)$. It follows from $|H:\Z(H)|=3^2$ that $\C_H(h)=\lg h\rg \Z(H)$ and $\C_H(h)$ is maximal in $H$.
Notice that $\forall t\in T-H$, $o(t)=3$. By Lemma \ref{gorder=3},
$[h^t,h]=1$. So $h^t\in \C_H(h)$. It follows that $\C_H(h)\trianglelefteq T$. Since $|T/ \C_H(h)|=3^2$, $T'\leq \C_H(h)$. Since $h$ is arbitrary, $T'\leq \Phi(H)=\Z(H)$.

Next, we claim $H\cong {\rm S}{\rm (}3^3{\rm )}$. Otherwise, $\exp(H)>\exp(\Z(H))$. Hence $a\not\in \Z(H)$ for all $a\in H$ with $o(a)=\exp(H)$.
Let $x\in T-H$. Then $o(x)=o(xa^{-1})=3$.
Since $T'\leq \Z(H)$, $$1=(xa^{-1})^3=x^3[x,a]^{{3}\choose{2}}[x,a,x][x,a,a]a^{-3}=[x,a]^3[x,a,x]a^{-3}.$$
Now, we consider $[x,a,x]$. Without loss of generality assume $H=\langle a,b\rangle$ and $H'=\lg c\rg$. Notice that
$[x,a]\in T'\leq \Z(H)=\Phi(H)=\lg a^3,b^3,c\rg$. We may assume that $[x,a]=a^{3i}b^{3j}c^k$.
Since $|H'|=3$, $H'\leq \Z(T)$. So $[c^k,x]=1$.
Hence $$[x,a,x]=[a^{3i}b^{3j}c^k,x]=[a^{3i}b^{3j},x]=[a^{3i},x][b^{3j},x]=[a^i,x]^3[b^j,x]^3.$$
Since $o(a)=\exp(H)>\exp(\Z(H))$,
$o([x,a]^3)<o(a^{-3})$ and $o([x,b]^3)<o(a^{-3})$.  So $o([x,a,x])<o(a^{-3})$.
It follow that $[x,a]^3[x,a,x]a^{-3}\neq1$. This is a contradiction. \qed %(又因为三者均属于$Z(M)$, 所以交换).
%So $H\cong {\rm S}{\rm (}3^3{\rm )}$. \qed

\begin{lemma}{\rm \cite[Theorem 4.5]{Zhang2019}}\label{p3feimetacyclic}
Assume $G$ is a finite nonabelian $p$-group and $p$ an odd prime.
Then all minimal nonabelian subgroups of $G$ are isomorphic to ${\rm S}{\rm (}p^3{\rm )}$  if and only if $\exp(G)=p$ or
$G=\H_p(G)\rtimes\langle a\rangle$, a semidirect product of
$\H_p(G)$ by $\langle a\rangle$, where the Hughes subgroup $\H_p(G)=B_1\times
B_2\times\cdots \times B_{n}$ is an abelian subgroup of index $p$
and $o(a)=p$. Moreover, for all $1\le i\le n$, $B_i\langle a\rangle$
is a group of maximal class with an abelian subgroup $B_i$ of index
$p$, and all elements of $B_i\langle a\rangle- B_i$ are of
order $p$.
\end{lemma}

\begin{lemma}\label{remark} Let $G$ be a nonabelian $p$-group of exponent $>p$ where $p\ge 5$. If all minimal nonabelian subgroups of $G$ are isomorphic to ${\rm S}{\rm (}p^3{\rm )}$, then $\N_G(S)$ is of exponent $p$ for any minimal nonabelian subgroup $S$ of $G$.
\end{lemma}
\demo
Assume the contrary and $\exp(\N_G(S))>p$ for some minimal nonabelian subgroup $S$ of $G$. Then there exists $a\in \N_G(S)- S$ such that $o(a)=p^2$. Let $L=S\lg a\rg$. Then $|L|\leq p^5$ and hence $\cl(L)\leq 4$. So $L$ is regular. By \cite[Proposition 10.28]{Ber1}, $L$ is generated by minimal nonabelian subgroups. So $\exp(L)=p$, a contradiction. \qed

\medskip
\begin{theorem}{\rm (A refinement of \cite[Theorem 4]{Janko})}\label{thm6}
 Let $G$ be a nonabelian $p$-group of exponent $>p$ which is not minimal nonabelian. Then all minimal nonabelian subgroups of $G$ are isolated in their normalizers if and only if all minimal nonabelian subgroups of $G$ are isomorphic to ${\rm S}{\rm (}p^3{\rm )}$, where $p\geqslant5$.
\end{theorem}
\demo $(\Leftarrow:)$ By Lemma \ref{remark}, $S$ is isolated in $\N_G(S)$ for any minimal nonabelian subgroup $S$ of $G$.

$(\Rightarrow:)$ By Lemma \ref{p>2}, $p>2$.
By Lemma \ref{njhqgl}, all minimal nonabelian subgroups of $G$ are isomorphic to ${\rm S}{\rm (}p^3{\rm )}$.
By Lemma \ref{p3feimetacyclic}, $G=\H_p(G)\rtimes\langle a\rangle$ has the structure described as Lemma \ref{p3feimetacyclic}.
Assume $p=3$. We will prove that this is a contradiction.

If $\cl (G)=2$, then $G$ is regular. Notice that all elements of $G- \H_3(G)$ are of order $3$ and $G=\lg G-\H_3(G)\rg$. We have $\exp(G)=3$.
This contradicts the hypothesis of $\exp(G)>3$.

If $\cl (G)\ge 3$, then there exists some $j$ such that $\cl(B_j\lg a\rg)\geq 3$.
Let $M=\Z_3(B_j\lg a\rg)\lg a\rg$, where $\Z_3(B_j\lg a\rg)$ is
the third term of upper central series of $B_j\lg a\rg$. Then $M$ is of maximal class, and $|M|=3^4$. By checking the list of groups of order $3^4$, there is a $c\in M$ such that $c^3\ne 1$. Let $S$ be a minimal nonabelian subgroup of $M$. Then $M\le \N_G(S)$. Since $c^3\cap S\neq 1$, $S$ is not isolated in its normalizer.  This is a contradiction.
\qed

\medskip

 \begin{theorem}{\rm(A refinement of \cite[Corollary 5]{Janko}\rm)}\label{Cro}
Let $G$ be a nonabelian $p$-group which is not minimal nonabelian. Then all minimal nonabelian subgroups of $G$ are isolated in $G$ if and only if $\exp(G)=p$.
\end{theorem}
\demo $(\Leftarrow:)$ It is obvious.

$(\Rightarrow:)$ Since all minimal nonabelian subgroups of $G$ are isolated in $G$, all minimal nonabelian subgroups of $G$ are isolated in their normalizers. By Lemma \ref{njhqgl}, all minimal nonabelian subgroups of $G$ are isomorphic to ${\rm S}{\rm (}p^3{\rm )}$.  % 下面来证明$p\geq5$.

Assume that $\exp(G)>p$. Then, by Lemma \ref{p3feimetacyclic}, $G=\H_p(G)\rtimes\langle a\rangle$ has the structure described as Lemma \ref{p3feimetacyclic}.
Since $\exp(G)>p$, there exists some $j$ such that $\exp(B_j\lg a\rg)>p$. So $\exp(B_j)>p$ and hence $\mho_1(B_j)\neq1$.
Since $B_j\lg a\rg$ is a group of maximal class, $|\Z(B_j\lg a\rg)|=p$. Hence $\Z(B_j\lg a\rg)\leq \mho_1(B_j)$.
Notice that $B_j$ is abelian. Then there exists $b\in B_j$ such that $\Z(B_j\lg a\rg)=\lg b^p\rg$.
Let $D=\Z_2(B_j\lg a\rg)\lg a\rg$, where $\Z_2(B_j\lg a\rg)$ is
the second term of upper central series of $B_j\lg a\rg$. Then $D$ is a nonabelian group of order $p^3$.
So $D\cong {\rm S}{\rm (}p^3{\rm )}$. Since $1\neq \lg b^p\rg=\Z(B_j\lg a\rg)\leq D$, $D$ is not isolated in $G$. This is a contradiction. \qed

\begin{rem}
{\rm (1)} If $G$ is the group in Theorem $\ref{thm6}$ and all minimal nonabelian subgroups of $G$ are isomorphic to ${\rm S}{\rm (}p^3{\rm )}$, then, by Lemma $\ref{p3feimetacyclic}$, $G$ has an abelian maximal subgroup $A$ of exponent $>p$ with $A ={\rm H}_p(G)$. By Lemma $\ref{remark}$, $\N_G(S)$ is of exponent $p$ for any minimal nonabelian subgroup $S$ of $G$. This means Theorem $\ref{thm6}$ contains all conclusions in {\rm\cite[Theorem 4(a)]{Janko}}.

{\rm (2)} The possibility {\rm (b)} in {\rm \cite[Theorem 4]{Janko}} does not occur. By Theorem $\ref{Cro}$, the $p$-groups satisfying
the hypothesis of {\rm \cite[Corollary 5]{Janko}} do not exist.
\end{rem}

\end {document}